\documentclass[11pt]{article}

\usepackage[a4paper,margin=2.5cm]{geometry}
\usepackage{amsmath,amssymb,amsfonts}
\usepackage{amsthm}
\usepackage{graphicx}
\usepackage{booktabs}
\usepackage{bm}
\usepackage{siunitx}
\usepackage{hyperref}
\usepackage{enumitem}
\usepackage{natbib}
\usepackage{setspace}
\usepackage{mathtools}
\usepackage{algorithm}
\usepackage{algpseudocode}

\usepackage{listings}
\usepackage{xcolor}

\lstdefinestyle{pythonstyle}{
  language=Python,
  basicstyle=\ttfamily\small,
  keywordstyle=\bfseries,
  commentstyle=\itshape\color{gray!70!black},
  stringstyle=\color{gray!50!black},
  showstringspaces=false,
  numbers=left,
  numberstyle=\tiny\color{gray},
  numbersep=8pt,
  frame=single,
  rulecolor=\color{black!30},
  framerule=0.4pt,
  breaklines=true,
  breakatwhitespace=true,
  tabsize=4,
  columns=flexible,
  xleftmargin=0pt,
  xrightmargin=0pt,
  belowskip=0.5em,
  aboveskip=0.5em,
}
\title{Dual-Network PINNs for Optimal Control: A Reproducible Benchmark on the Mass--Spring--Damper System}

\author{
Abdeladhim Tahimi \\
\small Centro de Engenharias e Ciências Agrárias, Universidade Federal de Alagoas, Brazil \\
\small\texttt{abdeladhim.tahimi@ceca.ufal.br}
\and
Rinaldo Vieira da Silva Junior \\
\small Centro de Engenharias e Ciências Agrárias, Universidade Federal de Alagoas, Brazil \\
\small\texttt{rinaldo.silva@ceca.ufal.br}
}

\date{\today}

\begin{document}

\maketitle

\begin{abstract}
This work presents a transparent and reproducible benchmark study of a direct
dual-network Physics-Informed Neural Network (PINN) formulation for the optimal
control of a mass--spring--damper system. The classical linear-quadratic optimal
control problem is solved by two independent classical methods -- Pontryagin's
Minimum Principle with single shooting, and direct transcription through
trapezoidal collocation -- and recast as a constrained optimization problem
solved by two feedforward neural networks: a state network whose boundary
conditions are enforced exactly through a composite cubic-and-mask ansatz, and
an unconstrained control network. The composite loss combines the physics
residual at the collocation points with a trapezoidal approximation of the
cost functional, weighted by a single scalar hyperparameter. On the benchmark
considered, the PINN reproduces the classical optimal cost to four significant
digits, satisfies the terminal state constraints exactly by construction, and
produces pointwise state and control errors that fall within the spread of the
two classical references. Training is approximately two orders of magnitude
slower than classical shooting on this benchmark, which is honestly reported.
The contribution is methodological clarity rather than methodological novelty:
the formulation and the accompanying Google Colab implementation are intended
to lower the barrier to entry for practitioners exploring PINN-based optimal
control without prior exposure to adjoint methods or two-point boundary value
problems.

\end{abstract}

\section{Introduction}
\label{sec:introduction}
% introduction.tex

Optimal control problems play a central role in mathematics, physics, and engineering because many physical systems are governed by differential equations while simultaneously requiring the optimization of constrained performance criteria. Applications arise in areas such as aerospace engineering, robotics, structural dynamics, energy systems, and fluid mechanics, where one seeks to determine control strategies capable of steering the evolution of a dynamical system toward a desired objective while respecting physical and operational constraints \citep{bryson1975applied,lions1971optimal}. Classical benchmark systems remain particularly valuable in this context because they provide analytically interpretable environments in which optimization methodologies can be studied rigorously. Among these benchmark problems, the mass--spring--damper system occupies a distinguished position due to its simplicity, physical interpretability, and broad relevance in vibration analysis and control theory.

Classical optimal control theory provides mathematically rigorous frameworks for addressing such constrained dynamical problems \citep{bryson1975applied,troltzsch2010optimal}. Direct methods reformulate the control problem into finite-dimensional optimization problems through discretization strategies, whereas indirect methods derive first-order optimality conditions based on variational principles and Pontryagin's Minimum Principle (PMP). These approaches have led to powerful computational techniques including shooting methods, collocation methods, pseudospectral methods, adjoint-based optimization, and Riccati formulations for linear-quadratic systems \citep{rao2009survey}. Nevertheless, despite their mathematical maturity and excellent performance for many benchmark problems, classical approaches may become increasingly difficult to implement in complex scenarios involving nonlinear dynamics, coupled multiphysics systems, inverse problems, sparse observations, or high-dimensional parameter spaces. In particular, the derivation and implementation of adjoint systems and two-point boundary value problems often represent a significant conceptual and computational barrier for newcomers entering the field.

In recent years, Physics-Informed Neural Networks (PINNs) have emerged as an important framework within scientific machine learning for solving differential-equation-based problems through physics-informed optimization. Introduced by \citet{raissi2019pinns}, PINNs approximate the solution of differential equations using neural networks trained through composite loss functions that incorporate physical residuals, boundary conditions, and observational constraints. The use of automatic differentiation enables the computation of derivatives directly through the computational graph, thereby allowing the governing equations to be embedded into the optimization process itself. Since their introduction, PINNs have been successfully applied to a wide variety of forward and inverse problems involving ordinary and partial differential equations, parameter estimation, hidden-physics discovery, and scientific data assimilation \citep{karniadakis2021piml,cuomo2022scientific}. Recent review papers have further highlighted the growing role of PINNs within the broader context of physics-informed scientific machine learning \citep{farea2024understanding,luo2025review,ren2025evolution}.

Beyond forward simulation problems, PINN methodologies are increasingly being extended toward optimal control and constrained optimization applications. Recent works have investigated direct and indirect PINN formulations for optimal control \citep{zhang2026directindirect}, differentiable optimization frameworks \citep{khimin2024autodiff}, Control PINNs \citep{barry2025controlpinn}, and physics-informed approaches for PDE-constrained optimization problems \citep{schiassi2024ponn}. These studies demonstrate that PINN formulations may naturally combine state approximation, control optimization, parameter estimation, and physical constraints within a unified differentiable optimization framework. Such formulations become particularly attractive in inverse problems and partially observed systems, where classical optimization pipelines may require repeated solver calls, mesh generation procedures, or iterative adjoint computations. Moreover, recent literature has emphasized the growing interaction between PINNs, differentiable programming, scientific machine learning, and control-oriented neural modeling \citep{wang2025survey,wang2026controlsurvey,legaard2023dynamical}.

Despite the rapid development of the field, the current literature has also become increasingly sophisticated from both mathematical and computational perspectives. Many recent contributions focus on high-dimensional PDE-constrained optimization, advanced neural architectures, operator-learning frameworks, Hamilton--Jacobi--Bellman formulations, or structure-preserving optimality systems \citep{barry2025controlpinn,zhang2026directindirect}. Although these developments are scientifically important, they may also create a substantial accessibility barrier for readers without prior experience in optimal control theory, adjoint methods, scientific machine learning, or differentiable programming. Recent educational and review-oriented contributions have attempted to improve the accessibility of PINN methodologies and automatic differentiation frameworks \citep{khimin2024autodiff,lawal2022review}. In particular, recent didactic work has emphasized the importance of exposing the internal mechanics of PINN training rather than treating automatic differentiation as a purely black-box procedure \citep{tahimi2026didactic,tahimi2026ad}. Nevertheless, transparent benchmark-oriented studies connecting classical optimal control formulations and direct PINN-based optimization strategies remain comparatively limited.

From an applied mathematics perspective, minimal and analytically verifiable benchmark problems remain essential for understanding the numerical and optimization behavior of emerging methodologies. Simple benchmark systems provide controlled environments in which approximation quality, optimization sensitivity, residual balancing, and training behavior can be studied transparently before extending the methodology toward more sophisticated applications. In this work, the classical mass--spring--damper system is formulated as a constrained optimal control problem with a quadratic objective functional. The benchmark is deliberately chosen because its mathematical structure is sufficiently simple to admit reliable classical reference solutions while remaining rich enough to illustrate the essential ingredients of direct PINN-based optimal control. The resulting framework allows a direct comparison between classical optimal-control formulations and a physics-informed optimization strategy employing separate neural representations for the state and control variables.

The objective of the present work is therefore not to propose a fundamentally new PINN architecture, nor to argue that PINNs outperform mature optimal-control solvers. Instead, this paper aims to provide a transparent, reproducible, and pedagogically structured benchmark study illustrating how constrained optimal-control problems can be reformulated within a direct physics-informed learning framework. More specifically, the work presents: (i) the mathematical formulation of a mass--spring--damper optimal control problem; (ii) its classical solution using standard optimal-control methodology; (iii) a direct dual-network PINN formulation employing independent neural representations for the state and control variables; (iv) a validation of the neural approximation against classical reference solutions; and (v) a comparative discussion of the advantages, limitations, and future potential of PINN-based optimal control formulations. In addition, all numerical experiments are designed to be reproducible through a companion Google Colab implementation intended to facilitate further educational and research-oriented exploration.

\section{Problem Formulation}
\label{sec:problem}
% problem_formulation.tex

The mass--spring--damper system is one of the most fundamental models in classical mechanics and vibration theory. It describes the motion of a point mass $m > 0$ attached to a linear spring of stiffness $k > 0$ and subject to viscous damping with coefficient $c > 0$, under the action of an external force $u(t)$ acting as the control input. The governing equation of motion is the second-order linear ordinary differential equation
\begin{equation}
    m\ddot{x}(t) + c\dot{x}(t) + k x(t) = u(t), \qquad t \in [0, T],
    \label{eq:msd}
\end{equation}
where $x(t) \in \mathbb{R}$ denotes the displacement of the mass from its equilibrium position, and dots denote differentiation with respect to time. The external forcing term $u(t) \in \mathbb{R}$ is treated as a time-varying control function to be determined through optimization.

The optimal control problem consists of finding a control function $u : [0,T] \to \mathbb{R}$ that minimizes a quadratic performance index measuring the integrated state deviation and control effort over the fixed time horizon $[0, T]$. Specifically, the objective functional is defined as
\begin{equation}
    J[u] = \int_0^T \left( x(t)^2 + r\, u(t)^2 \right) dt,
    \label{eq:cost}
\end{equation}
where $r > 0$ is a scalar weighting parameter balancing the relative penalization of state deviation and control effort. The control problem is subject to the dynamics~\eqref{eq:msd}, the initial conditions
\begin{equation}
    x(0) = x_0, \qquad \dot{x}(0) = v_0,
    \label{eq:ic}
\end{equation}
and the fixed terminal constraints
\begin{equation}
    x(T) = 0, \qquad \dot{x}(T) = 0,
    \label{eq:tc}
\end{equation}
requiring the system to reach a complete rest state at the final time $T$.

The complete optimal control problem is therefore stated as follows.

\begin{quote}
\textbf{Problem (P).} \textit{Find a control function $u \in L^2([0,T])$ such that the functional $J[u]$ defined in~\eqref{eq:cost} is minimized, subject to the dynamical constraint~\eqref{eq:msd}, the initial conditions~\eqref{eq:ic}, and the terminal constraints~\eqref{eq:tc}.}
\end{quote}

The numerical parameters adopted throughout this work are summarized in Table~\ref{tab:parameters}. The mass, spring constant, and damping coefficient are chosen as $m = k = c = 1$ in consistent SI units, yielding a damping ratio $\zeta = c / (2\sqrt{mk}) = 0.5$, which places the uncontrolled system in the underdamped regime. The time horizon is set to $T = 5\,\text{s}$, which is sufficiently long to observe both the natural transient response and the effect of the optimal control action. The initial conditions are $x_0 = 1\,\text{m}$ and $v_0 = 0\,\text{m/s}$, representing an initial displacement from equilibrium with zero velocity. The control weight is set to $r = 0.1$, which penalizes control effort lightly relative to state deviation, thereby encouraging active intervention by the optimizer. These parameter choices result in a well-posed problem admitting a unique optimal solution that serves as the reference throughout this study.

\begin{table}[ht]
\centering
\caption{Numerical parameters of the optimal control problem.}
\label{tab:parameters}
\begin{tabular}{llll}
\toprule
Parameter & Symbol & Value & Unit \\
\midrule
Mass                  & $m$ & $1$   & \si{\kilogram} \\
Spring constant       & $k$ & $1$   & \si{\newton\per\meter} \\
Damping coefficient   & $c$ & $1$   & \si{\newton\second\per\meter} \\
Time horizon          & $T$ & $5$   & \si{\second} \\
Initial displacement  & $x_0$ & $1$ & \si{\meter} \\
Initial velocity      & $v_0$ & $0$ & \si{\meter\per\second} \\
Terminal displacement & $x(T)$ & $0$ & \si{\meter} \\
Terminal velocity     & $\dot{x}(T)$ & $0$ & \si{\meter\per\second} \\
Control weight        & $r$ & $0.1$ & --- \\
\bottomrule
\end{tabular}
\end{table}

\section{Classical Optimal Control Solution}
\label{sec:classical}
% classical_solution.tex

\subsection{Optimality Conditions via Pontryagin's Minimum Principle}
\label{sec:pmp}

The application of Pontryagin's Minimum Principle requires recasting the second-order equation~\eqref{eq:msd} as an equivalent first-order system \citep{bryson1975applied}. Introducing the state variables $x_1(t) = x(t)$ and $x_2(t) = \dot{x}(t)$, the dynamics~\eqref{eq:msd} with $m = k = c = 1$ become
\begin{equation}
    \dot{x}_1(t) = x_2(t),
    \qquad
    \dot{x}_2(t) = -x_1(t) - x_2(t) + u(t),
    \label{eq:statespace}
\end{equation}
with initial conditions $x_1(0) = 1$, $x_2(0) = 0$, and fixed terminal constraints $x_1(T) = 0$, $x_2(T) = 0$.

To apply the Pontryagin Minimum Principle, one associates a costate variable $\lambda_i(t)$ with each state variable $x_i(t)$. The Hamiltonian of problem~\textbf{(P)} is defined as
\begin{equation}
    H(\bm{x}, u, \bm{\lambda}) =
    x_1^2 + r\,u^2
    + \lambda_1 x_2
    + \lambda_2 \left( -x_1 - x_2 + u \right),
    \label{eq:hamiltonian}
\end{equation}
where $\bm{x} = (x_1, x_2)^\top$ and $\bm{\lambda} = (\lambda_1, \lambda_2)^\top$.

\paragraph{Optimality condition.}
Since the Hamiltonian~\eqref{eq:hamiltonian} is strictly convex in $u$ and the control is unconstrained, the optimal control $u^*(t)$ is obtained by setting $\partial H / \partial u = 0$, which gives
\begin{equation}
    2r\,u^*(t) + \lambda_2(t) = 0
    \quad \Longrightarrow \quad
    u^*(t) = -\frac{\lambda_2(t)}{2r}.
    \label{eq:optimal_control}
\end{equation}

\paragraph{Costate equations.}
The costate dynamics are governed by
\begin{equation}
    \dot{\lambda}_1(t) = -\frac{\partial H}{\partial x_1} = -2x_1(t) + \lambda_2(t),
    \qquad
    \dot{\lambda}_2(t) = -\frac{\partial H}{\partial x_2} = -\lambda_1(t) + \lambda_2(t).
    \label{eq:costate}
\end{equation}
Since the terminal state is fully fixed by~\eqref{eq:tc}, no transversality conditions are required and the terminal costate values $\lambda_1(T)$ and $\lambda_2(T)$ are left free.

\paragraph{Two-point boundary value problem.}
Substituting the optimal control~\eqref{eq:optimal_control} into the state equations~\eqref{eq:statespace} and combining with the costate equations~\eqref{eq:costate} yields the following two-point boundary value problem (TPBVP): find $\bm{x}(t)$ and $\bm{\lambda}(t)$ satisfying
\begin{equation}
\frac{d}{dt}
\begin{pmatrix} x_1 \\ x_2 \\ \lambda_1 \\ \lambda_2 \end{pmatrix}
=
\begin{pmatrix}
    x_2 \\
    -x_1 - x_2 - \dfrac{\lambda_2}{2r} \\[6pt]
    -2x_1 + \lambda_2 \\
    -\lambda_1 + \lambda_2
\end{pmatrix},
\label{eq:tpbvp}
\end{equation}
subject to the boundary conditions
\begin{equation}
    x_1(0) = 1, \quad x_2(0) = 0,
    \qquad
    x_1(T) = 0, \quad x_2(T) = 0.
    \label{eq:tpbvp_bc}
\end{equation}
The system~\eqref{eq:tpbvp}--\eqref{eq:tpbvp_bc} constitutes a four-dimensional autonomous ODE with split boundary conditions: two conditions are prescribed at $t = 0$ and two at $t = T$. The unknown initial costate values $\lambda_1(0)$ and $\lambda_2(0)$ are determined by enforcing the terminal conditions numerically.

\paragraph{Numerical solution via shooting.}
The TPBVP is solved using a single shooting method \citep{rao2009survey}. Let $\bm{s} = (\lambda_1(0), \lambda_2(0))^\top$ denote the unknown initial costate vector. For a given $\bm{s}$, the system~\eqref{eq:tpbvp} is integrated forward from $t = 0$ to $t = T$ using a standard ODE solver, producing terminal state values $x_1(T;\bm{s})$ and $x_2(T;\bm{s})$. The shooting residual is defined as
\begin{equation}
    \bm{R}(\bm{s}) =
    \begin{pmatrix} x_1(T;\bm{s}) \\ x_2(T;\bm{s}) \end{pmatrix}
    -
    \begin{pmatrix} 0 \\ 0 \end{pmatrix},
    \label{eq:shooting_residual}
\end{equation}
and the correct initial costate is found by solving $\bm{R}(\bm{s}) = \bm{0}$ using a nonlinear root-finding algorithm. The resulting trajectories $x_1^*(t)$, $x_2^*(t)$, and $u^*(t) = -\lambda_2^*(t)/(2r)$ constitute the classical PMP reference solution to problem~\textbf{(P)}.

% -------------------------------------------------------------------
\subsection{Direct Transcription via Trapezoidal Collocation}
\label{sec:collocation}

Direct transcription methods reformulate the infinite-dimensional optimal control problem~\textbf{(P)} as a finite-dimensional nonlinear programming problem (NLP) by simultaneously discretizing the state and control trajectories on a fixed time grid \citep{rao2009survey}. Unlike the indirect approach of Section~\ref{sec:pmp}, no costate variables are introduced and no optimality conditions need to be derived analytically; instead, the physics constraints are enforced directly as algebraic equations at the collocation nodes.

Let $\{t_i\}_{i=0}^{N}$ denote a uniform partition of $[0, T]$ with step size $h = T/N$. Denote the discretized state and control values as $x_i \approx x(t_i)$, $v_i \approx \dot{x}(t_i)$, and $u_i \approx u(t_i)$. The state equations~\eqref{eq:statespace} are enforced at each interior interval $[t_i, t_{i+1}]$ using the trapezoidal rule:
\begin{align}
    x_{i+1} - x_i &= \frac{h}{2}\left( v_i + v_{i+1} \right),
    \label{eq:trap_x} \\
    v_{i+1} - v_i &= \frac{h}{2}\left[
        \left(-x_i - v_i + u_i\right)
        + \left(-x_{i+1} - v_{i+1} + u_{i+1}\right)
    \right],
    \label{eq:trap_v}
\end{align}
for $i = 0, \ldots, N-1$. The objective functional~\eqref{eq:cost} is approximated by the composite trapezoidal quadrature rule
\begin{equation}
    J_h = \frac{h}{2}\left[
        \left(x_0^2 + r\,u_0^2\right)
        + 2\sum_{i=1}^{N-1}\left(x_i^2 + r\,u_i^2\right)
        + \left(x_N^2 + r\,u_N^2\right)
    \right].
    \label{eq:trap_cost}
\end{equation}

The resulting NLP seeks the decision vector $\bm{z} = (x_0, v_0, \ldots, x_N, v_N, u_0, \ldots, u_N)^\top \in \mathbb{R}^{3(N+1)}$ that minimizes $J_h$ subject to the collocation constraints~\eqref{eq:trap_x}--\eqref{eq:trap_v}, the initial conditions $x_0 = 1$, $v_0 = 0$, and the terminal constraints $x_N = 0$, $v_N = 0$. This constrained NLP is solved using a standard gradient-based optimizer, with the collocation constraints enforced as equality constraints. The discrete trajectories $\{x_i^*\}$, $\{v_i^*\}$, and $\{u_i^*\}$ obtained at the solution provide a second independent reference against which both the PMP solution and the PINN approximation are subsequently validated.

\section{PINN Formulation}
\label{sec:pinn}
% pinn_formulation.tex

\subsection{Dual-Network Architecture}
\label{sec:architecture}

The PINN formulation of problem~\textbf{(P)} employs two independent feedforward neural networks: a state network $\mathcal{N}_x(t;\bm{\theta}_x)$ approximating the displacement trajectory $x(t)$, and a control network $\mathcal{N}_u(t;\bm{\theta}_u)$ approximating the control input $u(t)$. The use of separate networks for state and control variables is consistent with recent direct PINN formulations for optimal control \citep{barry2025controlpinn,zhang2026directindirect}. Both networks take the scalar time variable $t \in [0,T]$ as their sole input and produce a scalar output. Each network consists of $L$ hidden layers of width $W$, with hyperbolic tangent activations $\sigma(\cdot) = \tanh(\cdot)$ applied elementwise at each hidden layer. The output layer is linear. For a network with $L$ hidden layers, the raw output is defined recursively as
\begin{equation}
    \bm{h}^{(0)} = t, \qquad
    \bm{h}^{(\ell)} = \sigma\!\left( \bm{W}^{(\ell)} \bm{h}^{(\ell-1)} + \bm{b}^{(\ell)} \right),
    \quad \ell = 1, \ldots, L,
    \qquad
    \mathcal{F}(t;\bm{\theta}) = \bm{W}^{(L+1)} \bm{h}^{(L)} + \bm{b}^{(L+1)},
    \label{eq:network}
\end{equation}
where $\bm{\theta} = \{\bm{W}^{(\ell)}, \bm{b}^{(\ell)}\}$ collects all trainable weights and biases. The hyperbolic tangent activation is chosen because it is smooth, bounded, and admits derivatives of all orders, properties that are required for the automatic differentiation of the physics residual \citep{raissi2019pinns}. In the numerical experiments of Section~\ref{sec:results}, both networks are instantiated with $L = 4$ hidden layers and width $W = 32$, a configuration that is deliberately modest in order to preserve the didactic character of the benchmark.

\subsection{Hard Enforcement of Boundary Conditions}
\label{sec:ansatz}

The state network $\mathcal{N}_x$ must satisfy four boundary conditions exactly: the initial conditions $x(0) = 1$ and $\dot{x}(0) = 0$, and the terminal constraints $x(T) = 0$ and $\dot{x}(T) = 0$. Hard enforcement is achieved by constructing a composite ansatz of the form
\begin{equation}
    \hat{x}(t;\bm{\theta}_x) = \phi(t) + \psi(t)\,\mathcal{F}_x(t;\bm{\theta}_x),
    \label{eq:ansatz}
\end{equation}
where $\phi(t)$ is a fixed polynomial satisfying all four boundary conditions independently of the network output, $\psi(t)$ is a fixed masking function that vanishes at $t = 0$ and $t = T$ to first order, and $\mathcal{F}_x(t;\bm{\theta}_x)$ is the raw output of the state network as defined in~\eqref{eq:network}.

\paragraph{Construction of $\phi(t)$.}
A cubic polynomial $\phi(t) = a_0 + a_1 t + a_2 t^2 + a_3 t^3$ has four free coefficients, which are uniquely determined by imposing the four boundary conditions $\phi(0) = 1$, $\phi'(0) = 0$, $\phi(T) = 0$, $\phi'(T) = 0$. Solving the resulting linear system yields
\begin{equation}
    \phi(t) = 1 - \frac{3t^2}{T^2} + \frac{2t^3}{T^3}.
    \label{eq:phi}
\end{equation}
One may verify directly that $\phi(0) = 1$, $\phi'(0) = 0$, $\phi(T) = 0$, and $\phi'(T) = 0$, confirming that~\eqref{eq:phi} satisfies all required conditions exactly.

\paragraph{Construction of $\psi(t)$.}
The masking function must vanish at both endpoints to first order so that the composite ansatz~\eqref{eq:ansatz} preserves the boundary conditions of $\phi(t)$ regardless of the network output $\mathcal{F}_x$. The quadratic mask
\begin{equation}
    \psi(t) = t^2(T - t)^2 / T^4
    \label{eq:psi}
\end{equation}
satisfies $\psi(0) = \psi(T) = 0$ and $\psi'(0) = \psi'(T) = 0$, and is strictly positive on the open interval $(0, T)$, ensuring that $\mathcal{F}_x$ retains full expressive freedom in the interior of the time domain. The normalization factor $T^4$ keeps $\psi(t) \leq 1$ uniformly.

\paragraph{Verification.}
By construction, $\hat{x}(0;\bm{\theta}_x) = \phi(0) = 1$ and $\hat{x}(T;\bm{\theta}_x) = \phi(T) = 0$ for any value of $\bm{\theta}_x$. Differentiating~\eqref{eq:ansatz} and evaluating at the endpoints gives $\dot{\hat{x}}(0) = \phi'(0) + \psi'(0)\mathcal{F}_x(0) = 0$ and $\dot{\hat{x}}(T) = \phi'(T) + \psi'(T)\mathcal{F}_x(T) = 0$, confirming that all four boundary conditions are satisfied exactly and independently of the network parameters throughout training.

The control ansatz requires no hard enforcement. The control approximation is taken directly as the raw network output
\begin{equation}
    \hat{u}(t;\bm{\theta}_u) = \mathcal{F}_u(t;\bm{\theta}_u),
    \label{eq:control_ansatz}
\end{equation}
which is unconstrained and enters the loss function solely through the physics residual and the discretized cost functional defined in the following section.

\subsection{Loss Function}
\label{sec:loss}

The PINN loss function is constructed from two terms. The first term penalizes violation of the governing equation~\eqref{eq:msd} at a set of $N_c$ collocation points $\{t_i\}_{i=1}^{N_c}$ distributed uniformly on $[0, T]$. Using the composite ansatz~\eqref{eq:ansatz}, the second-order derivative $\ddot{\hat{x}}(t_i)$ is computed via automatic differentiation through the computational graph \citep{tahimi2026didactic}. The physics residual loss is
\begin{equation}
    \mathcal{L}_{\mathrm{phys}}(\bm{\theta}_x, \bm{\theta}_u) =
    \frac{1}{N_c} \sum_{i=1}^{N_c}
    \left[
        \ddot{\hat{x}}(t_i;\bm{\theta}_x)
        + \dot{\hat{x}}(t_i;\bm{\theta}_x)
        + \hat{x}(t_i;\bm{\theta}_x)
        - \hat{u}(t_i;\bm{\theta}_u)
    \right]^2.
    \label{eq:loss_phys}
\end{equation}
The second term approximates the objective functional~\eqref{eq:cost} using the composite trapezoidal rule on the same collocation grid:
\begin{equation}
    \mathcal{L}_{\mathrm{cost}}(\bm{\theta}_x, \bm{\theta}_u) =
    \frac{h}{2}\left[
        f(t_1) + 2\sum_{i=2}^{N_c - 1} f(t_i) + f(t_{N_c})
    \right],
    \qquad
    f(t_i) = \hat{x}(t_i;\bm{\theta}_x)^2 + r\,\hat{u}(t_i;\bm{\theta}_u)^2,
    \label{eq:loss_cost}
\end{equation}
where $h = T/(N_c - 1)$. The total loss is the weighted sum of both terms:
\begin{equation}
    \mathcal{L}(\bm{\theta}_x, \bm{\theta}_u) =
    \lambda_{\mathrm{phys}}\,
    \mathcal{L}_{\mathrm{phys}}(\bm{\theta}_x, \bm{\theta}_u)
    +
    \mathcal{L}_{\mathrm{cost}}(\bm{\theta}_x, \bm{\theta}_u),
    \label{eq:loss_total}
\end{equation}
where $\lambda_{\mathrm{phys}} > 0$ is a fixed scalar weight applied to the physics residual term. The introduction of $\lambda_{\mathrm{phys}}$ is motivated by the fundamentally different roles of the two loss terms. The cost term $\mathcal{L}_{\mathrm{cost}}$ is an approximation of the objective functional $J[u]$ that the optimizer seeks to minimize, while the physics residual $\mathcal{L}_{\mathrm{phys}}$ encodes the dynamical constraint~\eqref{eq:msd} and must be driven essentially to zero for the resulting trajectory to be a feasible solution of problem~\textbf{(P)}. Without a sufficiently large weight on the residual term, the optimizer is free to reduce $\mathcal{L}_{\mathrm{cost}}$ by violating the dynamics, producing a non-physical solution whose reported cost is artificially lower than the true optimum. The numerical value $\lambda_{\mathrm{phys}} = 10^{3}$ adopted in the present work is selected through a small number of preliminary training runs in which the relative magnitudes of $\lambda_{\mathrm{phys}}\,\mathcal{L}_{\mathrm{phys}}$ and $\mathcal{L}_{\mathrm{cost}}$ were monitored, and is held fixed for all experiments reported in Section~\ref{sec:results}. The role of $\lambda_{\mathrm{phys}}$ and the loss-balancing sensitivity that more generally characterizes PINN training are discussed further in Section~\ref{sec:discussion}.

\subsection{Training Procedure}
\label{sec:training}

Training consists of minimizing the total loss~\eqref{eq:loss_total} jointly with respect to $\bm{\theta}_x$ and $\bm{\theta}_u$ using a two-phase optimization strategy. In the first phase, the Adam optimizer is run for a fixed number of epochs with a scheduled learning rate, providing robust initial convergence from a random parameter initialization. In the second phase, the L-BFGS optimizer is applied to the parameters delivered by Adam, exploiting second-order curvature information to refine the solution toward a high-accuracy local minimum. This two-phase strategy is standard practice in the PINN literature \citep{raissi2019pinns,cuomo2022scientific} and is adopted here without modification. The collocation grid uses $N_c = 100$ uniformly spaced points on $[0, T]$. All derivatives required for~\eqref{eq:loss_phys} are computed via automatic differentiation using PyTorch's \texttt{autograd} engine, with \texttt{create\_graph=True} to enable higher-order differentiation through the ansatz. Network weights are initialized using the Glorot uniform scheme. The complete implementation is provided in the companion Colab notebook described in Section~\ref{sec:results}.

\section{Numerical Results and Validation}
\label{sec:results}
% numerical_results.tex

This section validates the dual-network PINN formulation of Section~\ref{sec:pinn} against the two independent classical references developed in Section~\ref{sec:pmp} and Section~\ref{sec:collocation}. All three methods solve the same optimal control problem (P) with the parameter values listed in Table~\ref{tab:parameters}. To enable a uniform comparison, the PMP shooting solution, the discrete collocation NLP solution (cubically interpolated from its $N = 100$ NLP grid), and the trained PINN are all evaluated on a common dense reference grid of $501$ uniformly spaced points on $[0, T]$. The composite trapezoidal rule on this dense grid is used to compute the cost functional and the $L^2$ error norms reported below. The full pipeline is implemented in a companion Google Colab notebook that mounts Google Drive, writes all numerical results to a single comma-separated values (CSV) file, and saves every figure through an independent rendering cell that re-reads the CSV from disk.

\subsection{Reproducibility and reference grid}
\label{sec:results_setup}

A single global random seed is used to initialize the random number generators of Python, NumPy, and PyTorch (including the CUDA backend when available). Network weights are initialized with the Glorot uniform scheme. Both PINN sub-networks employ four hidden layers of width $32$ and $\tanh$ activations, as specified in Section~\ref{sec:architecture}. The training collocation grid uses $N_c = 100$ points on $[0, T]$ and is held fixed throughout the optimization. All figures and statistical quantities reported below are reproduced bit-identically by the companion notebook provided the seed is preserved.

\subsection{Validation of the classical references}
\label{sec:results_classical}

The two classical methods provide independent reference solutions and are first validated against each other. The PMP shooting algorithm converges to terminal state residuals of order $10^{-15}$, and the collocation NLP enforces the terminal state as a hard equality constraint, yielding residuals of order $10^{-20}$. The two methods agree on the optimal cost to within $1.3 \times 10^{-3}$ in relative terms, with $J_{\mathrm{PMP}}^* = 6.872264 \times 10^{-1}$ and $J_{\mathrm{NLP}}^* = 6.881065 \times 10^{-1}$. The pointwise $L^\infty$ error between the two reference trajectories is $1.10 \times 10^{-3}$ for the displacement, $2.54 \times 10^{-3}$ for the velocity, and $1.10 \times 10^{-1}$ for the control. The small but non-negligible discrepancy reflects the trapezoidal discretization error of the NLP, which is of second order in the step size $h = T/N$; the PMP solution, integrated with a high-order adaptive Runge--Kutta scheme, is therefore treated as the gold reference in the remainder of this section.

\subsection{Training behaviour of the PINN}
\label{sec:results_training}

\begin{figure}[ht]
    \centering
    \includegraphics[width=0.85\textwidth]{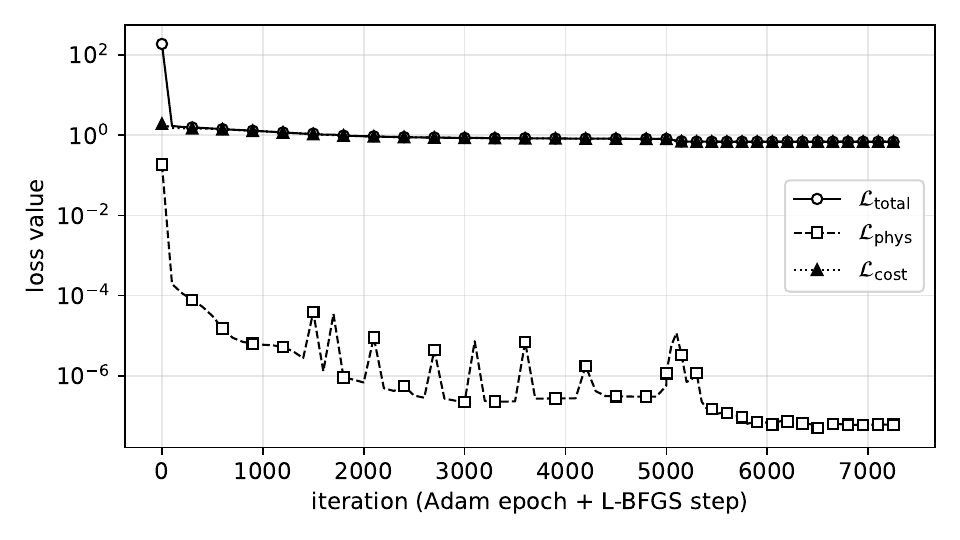}
    \caption{Training loss history of the dual-network PINN. The Adam phase occupies the first $5{,}000$ iterations; the L-BFGS refinement begins immediately afterwards.}
    \label{fig:loss_history}
\end{figure}

The PINN is trained by minimizing the composite loss $\mathcal{L} = \lambda_{\mathrm{phys}}\,\mathcal{L}_{\mathrm{phys}} + \mathcal{L}_{\mathrm{cost}}$ with $\lambda_{\mathrm{phys}} = 10^{3}$ following the discussion of Section~\ref{sec:loss}. A two-phase optimization strategy is used: $5{,}000$ Adam epochs with a constant learning rate of $10^{-3}$, followed by L-BFGS with strong Wolfe line search. The total training time on a standard Colab CPU runtime is approximately $125$ seconds. Figure~\ref{fig:loss_history} shows the evolution of the three loss components. The physics residual $\mathcal{L}_{\mathrm{phys}}$ drops by roughly six orders of magnitude during Adam and continues to decrease through the L-BFGS phase, reaching a final value of $6.17 \times 10^{-8}$. The cost term $\mathcal{L}_{\mathrm{cost}}$ converges to $6.874 \times 10^{-1}$ at the end of training, in close agreement with the PMP reference value. The total loss $\mathcal{L}$ is numerically indistinguishable from $\mathcal{L}_{\mathrm{cost}}$ once the residual has reached its asymptotic plateau, which is expected: $\lambda_{\mathrm{phys}} \mathcal{L}_{\mathrm{phys}} \approx 6 \times 10^{-5}$ is several orders of magnitude smaller than $\mathcal{L}_{\mathrm{cost}}$.

\subsection{Trajectory comparison}
\label{sec:results_trajectories}

\begin{figure}[ht]
    \centering
    \includegraphics[width=0.85\textwidth]{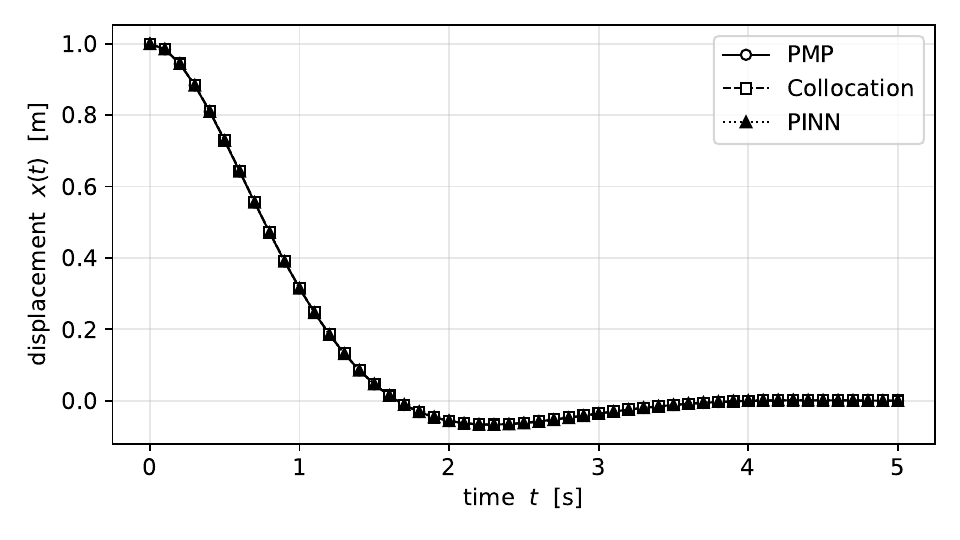}
    \caption{Optimal displacement $x(t)$ obtained by the three methods. The curves overlap on the plotting scale; the pointwise error is reported in Figure~\ref{fig:error_x}.}
    \label{fig:state_x}
\end{figure}

\begin{figure}[ht]
    \centering
    \includegraphics[width=0.85\textwidth]{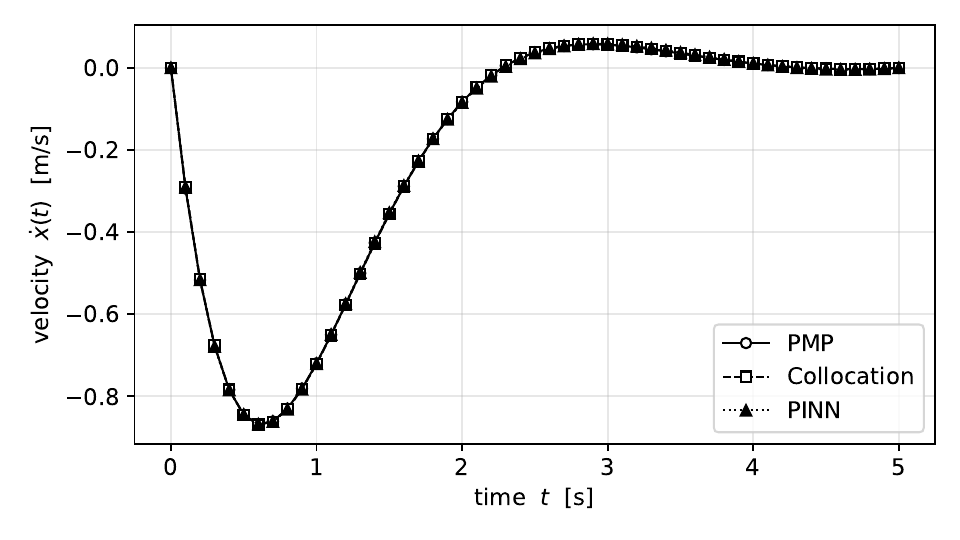}
    \caption{Optimal velocity $\dot{x}(t)$ obtained by the three methods. For the PINN, the velocity is obtained by automatic differentiation of the composite state ansatz $\hat{x}(t)$.}
    \label{fig:state_v}
\end{figure}

\begin{figure}[ht]
    \centering
    \includegraphics[width=0.85\textwidth]{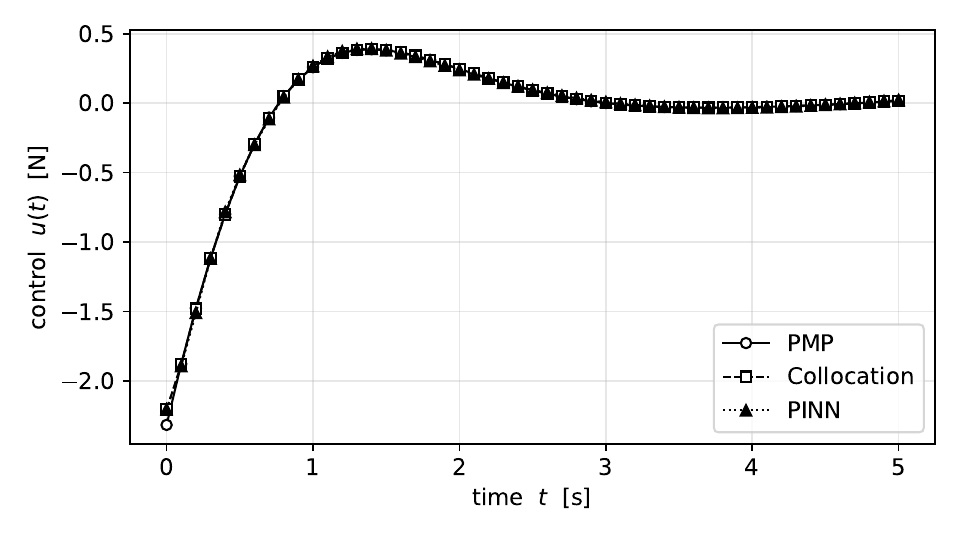}
    \caption{Optimal control $u(t)$ obtained by the three methods. A small visible discrepancy is observed near $t = 0$ where the magnitude of the optimal control is largest.}
    \label{fig:control_u}
\end{figure}

Figures~\ref{fig:state_x}, \ref{fig:state_v}, and \ref{fig:control_u} compare the three solution methods for the displacement, velocity, and control trajectories, respectively. The three curves are visually indistinguishable on the plotting scale of the state variables, and only the optimal control near the initial time exhibits a small visible discrepancy that is analysed in Section~\ref{sec:results_errors}. The optimal displacement decreases monotonically from $x_0 = 1\,\mathrm{m}$, briefly undershoots the equilibrium by approximately $0.07\,\mathrm{m}$ near $t = 2.5\,\mathrm{s}$, and returns asymptotically to zero. The optimal control begins at approximately $-2.3\,\mathrm{N}$, increases monotonically, crosses zero near $t = 1.1\,\mathrm{s}$, peaks near $+0.4\,\mathrm{N}$, and decays to zero at the final time. The PINN reproduces all of these qualitative features without any prior knowledge of the classical solutions.

\subsection{Pointwise error analysis}
\label{sec:results_errors}

\begin{figure}[ht]
    \centering
    \includegraphics[width=0.85\textwidth]{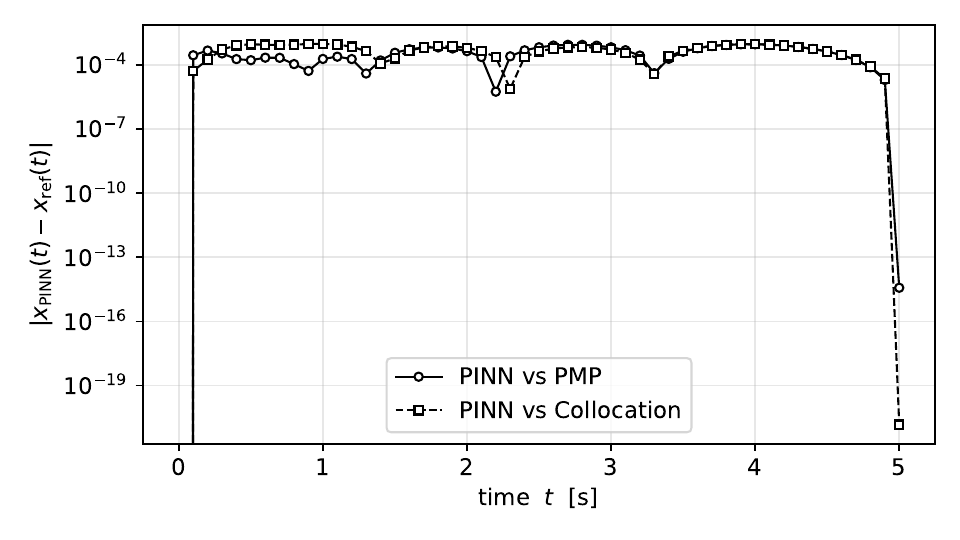}
    \caption{Absolute pointwise error of the PINN displacement, $\lvert \hat{x}(t) - x_{\mathrm{ref}}(t) \rvert$, against the PMP and collocation references. The collapse to machine precision at $t = 0$ and $t = T$ reflects the exact enforcement of the boundary conditions by the state ansatz.}
    \label{fig:error_x}
\end{figure}

Figure~\ref{fig:error_x} displays the absolute pointwise error of the PINN displacement against both classical references on a logarithmic vertical axis. Two features of the figure deserve emphasis. First, the interior error sits around $10^{-4}$--$10^{-3}\,\mathrm{m}$, which is comparable to the discrepancy between the two classical references themselves. The PINN therefore falls within the spread of the classical baseline. Second, the error collapses to machine precision at both endpoints. This is the visual signature of hard boundary-condition enforcement: the composite ansatz $\hat{x}(t) = \phi(t) + \psi(t)\,\mathcal{F}_x(t)$ guarantees $\hat{x}(0) = 1$ and $\hat{x}(T) = 0$ identically, regardless of the network parameters.

\begin{figure}[ht]
    \centering
    \includegraphics[width=0.85\textwidth]{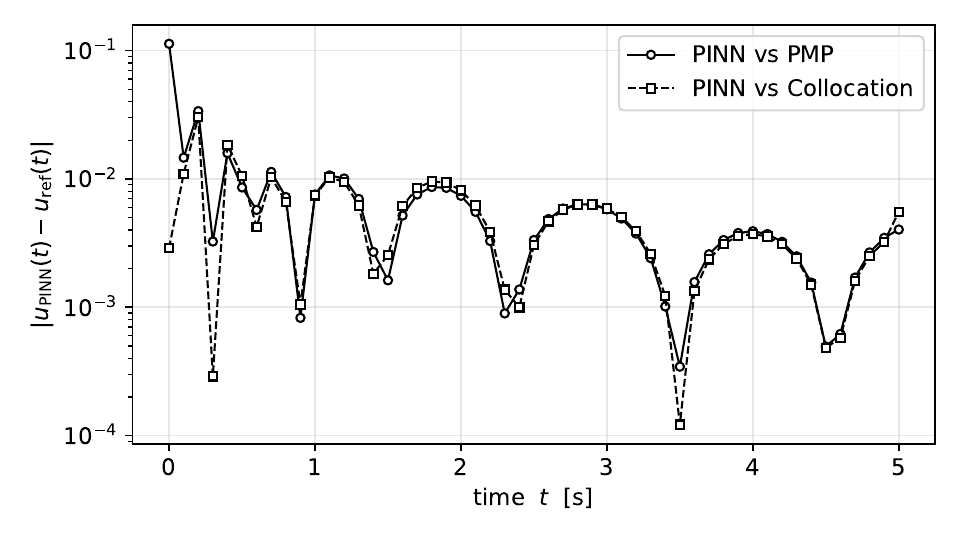}
    \caption{Absolute pointwise error of the PINN control, $\lvert \hat{u}(t) - u_{\mathrm{ref}}(t) \rvert$, against the PMP and collocation references. The control approximation is not subject to hard endpoint enforcement, and the error attains its maximum near $t = 0$ where the magnitude of the optimal control is largest.}
    \label{fig:error_u}
\end{figure}

Figure~\ref{fig:error_u} displays the corresponding error for the control. The error magnitude is approximately one order of magnitude larger than for the state, and is concentrated at the initial time. This behaviour has two contributing causes. First, the control approximation $\hat{u}(t)$ is the raw network output and is not constrained at the endpoints, in contrast to $\hat{x}$. Second, the magnitude of the optimal control is largest near $t = 0$, so a fixed relative error produces a larger absolute error there. Both effects are consistent with the hierarchy of approximation accuracy generally observed in PINN formulations, in which the directly approximated quantity is the most accurate, and quantities derived from it or independently approximated tend to exhibit progressively larger errors.

\subsection{Quantitative summary}
\label{sec:results_table}

Table~\ref{tab:results_costs} reports the optimal cost achieved by each method and the terminal constraint residuals. The PINN agrees with the PMP reference to four significant digits in the cost functional, and reaches the terminal state exactly by construction. Table~\ref{tab:results_errors} reports the $L^2$ and $L^\infty$ pointwise errors of all three pairwise comparisons. The errors of PINN against either classical reference are of the same order of magnitude as the errors between the two classical references themselves, confirming that the PINN solution lies within the spread of the classical baseline.

\begin{table}[ht]
\centering
\caption{Optimal cost and terminal constraint residuals for the three methods.}
\label{tab:results_costs}
\begin{tabular}{lccc}
\toprule
Method        & $J[u^*]$              & $\lvert x(T) \rvert$   & $\lvert \dot{x}(T) \rvert$ \\
\midrule
PMP           & $6.872264 \times 10^{-1}$ & $3.76 \times 10^{-15}$ & $9.09 \times 10^{-15}$ \\
Collocation   & $6.881065 \times 10^{-1}$ & $1.48 \times 10^{-21}$ & $8.47 \times 10^{-20}$ \\
PINN          & $6.871579 \times 10^{-1}$ & $0$                    & $0$                    \\
\bottomrule
\end{tabular}
\end{table}

\begin{table}[ht]
\centering
\caption{Pointwise $L^2$ and $L^\infty$ errors of the three pairwise comparisons, computed on the dense reference grid of $501$ uniformly spaced points.}
\label{tab:results_errors}
\begin{tabular}{llcc}
\toprule
Comparison              & Quantity      & $L^2$ error            & $L^\infty$ error       \\
\midrule
PINN vs PMP             & $x(t)$        & $1.14 \times 10^{-3}$  & $9.29 \times 10^{-4}$  \\
                        & $\dot{x}(t)$  & $3.37 \times 10^{-3}$  & $3.70 \times 10^{-3}$  \\
                        & $u(t)$        & $2.43 \times 10^{-2}$  & $1.13 \times 10^{-1}$  \\
\midrule
PINN vs Collocation     & $x(t)$        & $1.36 \times 10^{-3}$  & $9.74 \times 10^{-4}$  \\
                        & $\dot{x}(t)$  & $3.51 \times 10^{-3}$  & $3.23 \times 10^{-3}$  \\
                        & $u(t)$        & $1.91 \times 10^{-2}$  & $4.15 \times 10^{-2}$  \\
\midrule
Collocation vs PMP      & $x(t)$        & $9.88 \times 10^{-4}$  & $1.10 \times 10^{-3}$  \\
                        & $\dot{x}(t)$  & $1.59 \times 10^{-3}$  & $2.54 \times 10^{-3}$  \\
                        & $u(t)$        & $1.19 \times 10^{-2}$  & $1.10 \times 10^{-1}$  \\
\bottomrule
\end{tabular}
\end{table}

Wall-clock training time is reported for completeness. On the Colab CPU runtime, PMP shooting required $1.2\,\mathrm{s}$, the SLSQP collocation NLP required $10.4\,\mathrm{s}$, and PINN training required $125.0\,\mathrm{s}$ in total ($79.2\,\mathrm{s}$ for the Adam phase and $45.8\,\mathrm{s}$ for the L-BFGS phase). The PINN is therefore not competitive with the classical methods in training time on this benchmark, a point discussed further in Section~\ref{sec:discussion}.

\subsection{Selected code excerpts}
\label{sec:results_listings}

For pedagogical clarity, three minimal Python excerpts from the companion notebook are reproduced here. The complete implementation is available in the Colab notebook.

The composite ansatz of Section~\ref{sec:ansatz} is implemented in three lines:

\begin{lstlisting}[caption={Hard enforcement of the four boundary conditions via the composite ansatz.},label={lst:ansatz}]
def phi_fn(t):
    return 1.0 - 3.0 * (t / T) ** 2 + 2.0 * (t / T) ** 3

def psi_fn(t):
    return (t ** 2) * (T - t) ** 2 / (T ** 4)

def x_hat(net_x, t):
    return phi_fn(t) + psi_fn(t) * net_x(t)
\end{lstlisting}

The composite loss combines the physics residual at the collocation points with the trapezoidal quadrature of the cost functional. The required derivatives of $\hat{x}(t)$ are obtained by automatic differentiation through the ansatz with \verb|create_graph=True| to enable backpropagation of second-order derivatives:

\begin{lstlisting}[caption={Composite loss function combining physics residual and trapezoidal cost.},label={lst:loss}]
def losses():
    x_val = x_hat(net_x, t_col_pinn)
    u_val = u_hat(net_u, t_col_pinn)
    dx  = torch.autograd.grad(x_val, t_col_pinn,
              torch.ones_like(x_val), create_graph=True)[0]
    ddx = torch.autograd.grad(dx, t_col_pinn,
              torch.ones_like(dx), create_graph=True)[0]
    residual = ddx + dx + x_val - u_val
    L_phys = (residual ** 2).mean()
    f = x_val ** 2 + r * u_val ** 2
    L_cost = 0.5 * h_pinn * (f[0] + 2.0 * f[1:-1].sum() + f[-1]).squeeze()
    return L_phys, L_cost
\end{lstlisting}

The optimization loop combines an Adam warm-up phase with an L-BFGS refinement phase, following standard PINN practice:

\begin{lstlisting}[caption={Two-phase optimization: Adam warm-up followed by L-BFGS refinement.},label={lst:training}]
# Phase 1: Adam
for epoch in range(ADAM_EPOCHS):
    opt_adam.zero_grad()
    L_phys, L_cost = losses()
    L = LAMBDA_PHYS * L_phys + L_cost
    L.backward()
    opt_adam.step()

# Phase 2: L-BFGS
def closure():
    opt_lbfgs.zero_grad()
    L_phys, L_cost = losses()
    L = LAMBDA_PHYS * L_phys + L_cost
    L.backward()
    return L
opt_lbfgs.step(closure)
\end{lstlisting}

\section{Discussion}
\label{sec:discussion}
% discussion.tex

The results of Section~\ref{sec:results} establish that the dual-network PINN formulation reproduces the classical optimal control solution of problem~\textbf{(P)} to four significant digits in the cost functional and to three significant digits in the pointwise state trajectory, on a benchmark for which two independent classical references are available. The purpose of the present section is to interpret these results in their proper methodological context, to acknowledge the limitations of the formulation studied here, and to identify the regimes in which a direct PINN approach to optimal control may become genuinely competitive with classical methodologies.

\subsection{Loss balancing and the role of the residual weight}
\label{sec:disc_weighting}

The introduction of the weighting parameter $\lambda_{\mathrm{phys}}$ in Equation~\eqref{eq:loss_total} is the single most consequential design choice in the formulation. Preliminary experiments without weighting confirmed the failure mode anticipated in Section~\ref{sec:loss}: the optimizer reduced the cost term $\mathcal{L}_{\mathrm{cost}}$ to values substantially below the true optimum $J^{*}$, at the price of a physics residual $\mathcal{L}_{\mathrm{phys}}$ of order $10^{-3}$ that produced a non-physical trajectory. The resulting solution is inadmissible: it minimizes a quadratic functional over a function space in which the dynamical constraint is only weakly satisfied. The introduction of a fixed weight $\lambda_{\mathrm{phys}} = 10^{3}$ restored the correct hierarchy between constraint satisfaction and objective minimization, driving the residual below $10^{-7}$ while leaving the cost term in close agreement with both classical references.

This behaviour is consistent with the broader experience reported in the PINN literature, in which loss-balancing has been identified as a recurrent source of training pathologies \citep{cuomo2022scientific,farea2024understanding}. Several systematic approaches have been proposed in recent years, including adaptive weighting strategies, neural tangent kernel diagnostics, and constraint-based reformulations through Lagrangian or penalty-method updates \citep{ren2025evolution,luo2025review}. The present work deliberately employs the simplest possible weighting strategy -- a single fixed scalar -- in order to preserve the didactic character of the benchmark and to make the role of the residual weight transparent. More sophisticated weighting schemes are likely to be necessary in extensions to nonlinear dynamics or higher-dimensional state spaces, where the asymmetry between the two loss components may evolve substantially during training.

\subsection{Hard versus soft boundary condition enforcement}
\label{sec:disc_ansatz}

The composite ansatz $\hat{x}(t) = \phi(t) + \psi(t)\,\mathcal{F}_x(t;\bm{\theta}_x)$ enforces all four boundary conditions of problem~\textbf{(P)} exactly, independently of the network parameters. The benefit of this construction is visible in Figure~\ref{fig:error_x}: the pointwise displacement error collapses to machine precision at $t = 0$ and $t = T$, and the terminal constraints reported in Table~\ref{tab:results_costs} are satisfied identically. In contrast, the control approximation $\hat{u}(t)$ is the raw network output and is not constrained at the endpoints. The largest absolute error in the control accordingly occurs at $t = 0$, where the magnitude of the optimal control is greatest and the approximation enjoys no structural support from the ansatz.

This asymmetric design choice is intentional. In the problem considered here, the state is subject to explicit boundary conditions, whereas the control admits no natural endpoint constraint. Extending the hard-enforcement strategy to the control would require either the introduction of artificial endpoint conditions that are not present in the original problem formulation, or the construction of problem-specific ansatz functions whose role and motivation are less easily generalized. The trade-off between accuracy at the endpoints and faithfulness to the original problem formulation is itself a useful pedagogical observation, and extending the present framework to richer classes of constraints -- such as inequality bounds on the control, or terminal cost functionals with no terminal state constraint -- remains an interesting direction for further investigation.

\subsection{Approximation hierarchy and the cost of derived quantities}
\label{sec:disc_hierarchy}

The pointwise errors of Table~\ref{tab:results_errors} exhibit a clear hierarchy in which the directly approximated displacement is the most accurate, the velocity obtained by automatic differentiation through the ansatz is approximately one order of magnitude less accurate, and the independently approximated control is yet another order of magnitude less accurate. This hierarchy is a consequence of two distinct mechanisms. First, derivatives of a neural approximation $\hat{x}(t)$ are typically less accurate than the approximation $\hat{x}(t)$ itself, since the training loss applies pressure on the approximation but not directly on its derivatives. Second, the control network operates without the structural support of the ansatz, and is therefore the least constrained of the two networks in problem~\textbf{(P)}.

This hierarchy has practical implications for problem formulations in which the quantities of primary interest are not the directly approximated ones. In tracking problems, for example, the cost functional may depend explicitly on derivatives of the state, in which case the formulation studied here would impose an accuracy penalty on the most critical quantity. Reformulations in which both the state and its derivatives are independently approximated -- at the cost of additional networks and additional loss terms -- have been investigated in the recent PINN literature and may be appropriate in such cases \citep{barry2025controlpinn}.

\subsection{Computational cost and the question of competitiveness}
\label{sec:disc_cost}

The wall-clock training time reported in Section~\ref{sec:results_table} indicates that the PINN is two orders of magnitude slower than PMP shooting and one order of magnitude slower than the trapezoidal collocation NLP on this benchmark. This finding is not, on its own, an argument against PINN methods. It is rather an acknowledgment that the benchmark is a setting in which classical methods are already at their most efficient: a linear-quadratic problem of low dimension and short horizon, with a smooth unique optimum and exactly two terminal constraints \citep{bryson1975applied,rao2009survey}. Under such conditions, neither shooting nor collocation faces any of the difficulties for which PINN formulations were originally introduced.

The genuine question is not whether PINNs outperform classical methods on problems for which classical methods are already optimal, but whether they remain applicable in regimes where classical methods become difficult or impractical to deploy. Several such regimes have been identified in the literature. PINN formulations naturally accommodate sparse or noisy observations through additional data-fidelity loss terms \citep{raissi2019pinns,karniadakis2021piml}, do not require the derivation of adjoint equations for inverse or parameter-identification problems \citep{khimin2024autodiff}, scale gracefully to high-dimensional state spaces in which traditional collocation becomes intractable, and integrate seamlessly with end-to-end differentiable pipelines \citep{wang2025survey,wang2026controlsurvey}. None of these advantages is exercised by the benchmark of the present work, which was selected to maximize interpretability rather than to showcase the strengths of the PINN paradigm.

\subsection{Limitations and scope of the present study}
\label{sec:disc_limitations}

Several limitations should be acknowledged. The benchmark is linear, scalar, finite-horizon, and admits a smooth and unique optimum; none of these properties holds in many applications of practical interest. The weighting parameter $\lambda_{\mathrm{phys}}$ is fixed by inspection rather than chosen adaptively, and the training procedure relies on a standard two-phase Adam plus L-BFGS strategy without architectural innovations such as residual blocks, Fourier features, or domain decomposition. A single random seed is used to produce the results reported here; while the qualitative behaviour is robust across seeds, no systematic study of training variability is undertaken. Finally, the collocation grid of the PINN, the collocation grid of the NLP, and the dense reference grid used for error evaluation are all uniform, and no investigation of non-uniform or adaptive grid placement is reported. Each of these limitations defines a natural direction for further investigation.

The contribution of this work is therefore best understood as a transparent and reproducible benchmark study rather than a methodological advance. Its value lies in establishing a clear bridge between the classical optimal control literature and the rapidly growing literature on direct PINN formulations for constrained optimization \citep{barry2025controlpinn,zhang2026directindirect,schiassi2024ponn}, in a setting where every numerical claim can be independently verified against an analytically interpretable reference. The companion Colab notebook is intended to lower the barrier to entry for practitioners who wish to explore PINN-based optimal control without first mastering the full machinery of adjoint methods and two-point boundary value problems.

\section{Conclusion}
\label{sec:conclusion}
% conclusion.tex

This work has presented a transparent and reproducible benchmark study of a
direct dual-network Physics-Informed Neural Network formulation for the optimal
control of a mass--spring--damper system. The classical optimal control problem
was formulated explicitly, solved by two independent classical methods --
Pontryagin's Minimum Principle with single shooting, and direct transcription
through trapezoidal collocation -- and reformulated as a constrained optimization
problem to be solved by a pair of feedforward neural networks. Hard boundary
condition enforcement on the state was achieved through a composite ansatz built
from a cubic polynomial and a quadratic mask, while the control was approximated
by an unconstrained network. The composite loss combined a physics residual at
the collocation points with a trapezoidal approximation of the cost functional,
weighted by a single scalar hyperparameter selected by direct inspection.

On the benchmark studied here, the PINN approximation reproduces the classical
optimal cost to four significant digits, satisfies the terminal state constraints
exactly by construction, and produces pointwise state and control errors that
fall within the spread of the two classical references. Training time on a
standard Colab CPU runtime is approximately two minutes and is two orders of
magnitude larger than the cost of the classical methods on the same problem.
These findings establish that direct PINN formulations can reach
classical-grade accuracy on a benchmark optimal control problem without any
modification of the basic physics-informed training framework, and without the
explicit derivation of adjoint equations or two-point boundary value problems.

The contribution of the present work is methodological clarity rather than
methodological novelty. The benchmark was deliberately selected to admit
verifiable classical references and to expose the essential ingredients of a
direct PINN-based optimal control formulation in a setting where each design
choice can be examined in isolation. The companion Google Colab notebook, all
numerical results, and all figures are made available to support reproduction,
extension, and educational use of the framework. Natural directions for further
investigation include nonlinear and higher-dimensional dynamics, inequality
constraints on the control, adaptive loss weighting and collocation strategies,
and problem formulations in which the classical adjoint approach becomes
genuinely difficult to deploy. In such regimes, the qualitative advantages of
direct PINN-based optimization -- avoidance of adjoint derivations, natural
handling of partial observations, and seamless integration with differentiable
pipelines -- are expected to translate into practical competitiveness with
classical methodologies.

\section*{Code and Data Availability}
\label{sec:availability}
% code_data_availability.tex

This work does not involve experimental or observational datasets. All results are generated computationally from the differential equation~\eqref{eq:msd}, the boundary conditions of problem~\textbf{(P)}, and the parameter values of Table~\ref{tab:parameters}. To guarantee full reproducibility, the following materials are provided.

\paragraph{Source code.}
The complete Python/PyTorch implementation is provided as a single Jupyter notebook designed for Google Colab. The notebook mounts Google Drive, creates the output folders, implements the three solution methods (PMP shooting, trapezoidal collocation, dual-network PINN), trains the PINN with a fixed random seed, evaluates all three methods on a common dense reference grid, and produces every figure and table appearing in Section~\ref{sec:results}. Each figure cell is self-contained and reads the saved numerical results from disk, so any single figure can be regenerated without re-executing the training pipeline.

\paragraph{Computational artifacts.}
All numerical results reported in Section~\ref{sec:results} are written to two comma-separated values (CSV) files: \texttt{trajectories.csv} contains the state, velocity, control, and pairwise pointwise errors of the three methods on the dense reference grid of $501$ uniformly spaced points; and \texttt{loss\_history.csv} records the evolution of the three loss components during the Adam and L-BFGS training phases. Every numerical value reported in the manuscript tables can be reproduced from these two files alone.

\paragraph{Archive.}
The complete codebase, including the notebook and all computational artifacts, is publicly archived with a permanent DOI via Zenodo~\citep{tahimi2026notebook}.

\section*{Acknowledgments}
\label{sec:acknowledgments}
% acknowledgments.tex

During the development of this work, the authors used AI-based language models (Claude, Anthropic, 2024--2026) as assistive tools for writing, code development, and computational verification. All conceptual and scientific decisions --- including the choice of benchmark problem, the design of the dual-network PINN formulation, the selection of validation strategy against two independent classical references, and the critical discussion of the method's scope and limitations --- originated with the authors. All mathematical derivations, numerical values, and computational results were independently verified by the authors, who assume full responsibility for the integrity, accuracy, and originality of the submitted work.

% -------------------------------------------------
% Bibliography
% -------------------------------------------------
\bibliographystyle{plainnat}
\bibliography{references}

\end{document}